\documentclass{drna}

\usepackage{latexsym,todonotes}
\usepackage{dirtytalk}
\usepackage{xcolor}

\usepackage{graphicx}

\usepackage{caption,subcaption}
\usepackage{geometry}
\usepackage{float}

\usepackage{mathtools}

\usepackage{cellspace} 
\usepackage[style=numeric, sorting=none, url=false, doi=false, eprint=false, isbn=false]{biblatex} %Imports biblatex package
\begin{document}

\title{Fast-Decaying Polynomial Reproduction}

\author{Cappellazzo}{Giacomo}{a}
\author{De Marchi}{Stefano}{b}

\affiliation{a}{Department of Mathematics ``Tullio Levi-Civita", University of Padova, Italy}
\affiliation{b}{Department of Medicine (DIMED), University of Padova, Italy}

% The following two commands are filled by the publisher
\volume{17}
\firstpage{1}

\maketitle

%%%%%%%%%%%%%%%%%%%%%%%%%%%%%%%%%%%%%%%%%%%%%%%%%%%%%%%%%%%%%%%%%%%%

\begin{abstract}
Polynomial reproduction plays a crucial role in deriving error estimates for various approximation schemes. In particular, local polynomial reproduction is a key ingredient in both error estimation and stability analysis. However, for certain computationally relevant methods, such as Rescaled Localized Radial Basis Functions (RL-RBF), this requirement constitutes a limitation.
To enable the analysis of a broader class of approximation methods in a unified and efficient manner, the present work introduces a framework based on fast-decaying polynomial reproduction. In this approach, we do not restrict ourselves to compactly supported basis functions. Instead, we allow the basis functions to decay to zero at infinity, with the decay rate controlled as a function of the \textit{separation distance}. The adoption of fast-decaying polynomial reproduction yields stable and convergent approximation schemes. These methods can achieve smoothness when used in conjunction with moving least squares.
All theoretical results presented in this paper regarding the rate of convergence, the Lebesgue constant and the smoothness of the approximant have been numerically validated, including in the multivariate setting.
\end{abstract}

%%%%%%%%%%%%%%%%%%%%%%%%%%%%%%%%%%%%%%%%%%%%%%%%%%%%%%%%%%%%%%%%%%%%

\section{Introduction}

Polynomials are widely employed in approximation theory owing to their straightforward definition, ease of implementation and the availability of direct error estimates through Taylor expansions. Local polynomial reproduction plays a fundamental role in error estimation for various approximation schemes (cf. \cite{schaback_1997, Flyer2006}). In particular, the author in \cite{wendland_2001} derived a uniform bound on the Lebesgue constants for compactly supported radial basis functions, which is independent of the space dimension. Nevertheless, oversampling may still occur: the number of nodes in the support of each basis function must satisfy the constraints imposed by the dimension of the polynomial space being reproduced; as a result, the cardinality of the nodes can increase significantly compared with classical polynomial reproduction.
\\
The insights from arguments presented in \cite{demarchi2020} allow us to reformulate the notion of polynomial reproduction within a more general framework. This perspective enables a more natural and effective contextualization of the Rescaled Localized Radial Basis Function (RL-RBF) method, originally introduced in \cite{deparis_2014}, through the concept of fast-decaying polynomial reproduction. While the convergence of the RL-RBF interpolant was analyzed in \cite{demarchi2020}, that work left open a conjecture concerning the sum of the cardinal functions associated with the RL-RBF approximant.
\\
\\
This paper introduces a new approach to polynomial reproduction based on the assumption that the basis functions decay to zero at infinity. For a comprehensive overview of the field, we also refer to the recent monograph \cite{Buhmann2022}, which provides a detailed account of the latest developments in quasi-interpolation techniques. The flexibility of the fast-decaying polynomial reproduction framework further allows for an alternative perspective on these schemes. Quasi-interpolation methods using functions with non-compact support have been thoroughly analyzed in the literature, yielding a precise characterization of the approximation spaces for functions in the Sobolev space $W^{k}_{\infty}$ (cf. \cite{LightCheney1992}). Under suitable conditions (specifically, algebraic decay combined with possible singularity at the origin) our framework recovers the construction presented in \cite{LightCheney1992}.
\\
In the present work, the Gaussian kernel plays a central role, since its behavior at infinity can be precisely \textit{controlled}. This idea resonates with the concept of \textit{controlled approximation} previously explored by G. Strang \cite{Strang70}. Explicit error expansions for \textit{approximate approximations} using Gaussian kernels are provided in \cite{MazyaSchmidt1996}. Although this approach allows the construction of highly accurate approximations, the resulting approximants generally do not converge. The authors showed that convergence can only be achieved when the method parameters depend on the distribution of the data points in the domain. It is within this context that the present work introduces improvements.
A central role in this work is played by the Moving Least Squares technique, whose versatility also allows for the approximation of functions with discontinuities.  \cite{demarchi_2023a} proposes a moving least-squares scheme in which the weights are constructed from variably scaled discontinuous kernels. Furthermore, a coordinate-free moving least squares method that acts directly on point clouds, without tangent-plane approximations for algebraic manifolds, is presented in \cite{Hangelbroek2026}.
\\
The paper is organized as follows. In Section \ref{sec_preliminaries}, we recall the necessary definitions and introduce the concept of local polynomial reproduction in quasi-uniform settings, together with several auxiliary results that will be used in the subsequent discussion. In Section \ref{sec_stability_convergence}, we present the motivation for the proposed fast-decaying polynomial reproduction framework and establish its convergence and stability properties. In Section \ref{section_ex_c_infinity_fast}, we apply this framework to the \textit{moving least squares method}. 
In Section \ref{sec_numerical_experiments}, we present and discuss several numerical experiments. These tests compare the performance of the proposed approximation methods in both univariate and multivariate settings.
The analysis is completed by examining the behavior of basis functions with algebraic decay, which exhibit a singularity at the origin and thus represent a degenerate case for the proposed framework. Finally, in Section 7, we summarize the main results of the paper and outline possible directions for future research.

\section{Notations and preliminaries}
\label{sec_preliminaries}
We begin by recalling the notion of {\it local polynomial reproduction} (cf., \cite[Definition 3.1]{wendland_2004}). Let $X=\{x_1,\dots,x_N\} \subseteq \Omega\subseteq \mathbb{R}^d$ be a finite set of distinct points. Their {\it separation radius} and {\it fill distance} are defined as
\begin{equation*}
    q_{X} \coloneqq \frac{1}{2} \min_{i \neq j} \|x_i-x_j\|_2, \qquad h_{X,\Omega} \coloneqq \sup_{x \in \Omega} \min_{1 \leq j \leq N} \|x-x_j\|_{2}.
\end{equation*}
\begin{definition}\label{def_local_app_procc}
A numerical scheme that produces for every set $X=\{ x_1, \dots, x_N\} \subseteq \Omega$ a family of functions $u_j = u_j^X : \Omega \rightarrow \mathbb{R}, \, 1 \leq j \leq N,$ provides local polynomial reproduction of degree $\ell$ on $\Omega$ if there exist constants $h_0, C_1, C_2 \in \mathbb{R}_{>0}$ such that
	\begin{enumerate}
		\item $\displaystyle \sum_{j=1}^{N} p(x_j)u_j = p$ for each $p \in \pi_{\ell}(\mathbb{R}^d)$,
		\item $\displaystyle \sum_{j=1}^{N} |u_j(x)| \leq C_1$ for all $x \in \Omega$,
		\item $u_j(x)=0$ if $\|x-x_j\|_2 > C_2 h_{X,\Omega}$
	\end{enumerate}
are satisfied for all $X$ with $h_{X,\Omega} \leq h_0$. 
\end{definition}

\noindent
Note that  $\pi_{\ell}(\mathbb{R}^d)$ denotes the space of polynomials of degree at most $\ell$ in $\mathbb{R}^d$.

\vskip 0.1in
\noindent Given a function $f: \Omega \rightarrow \mathbb{R}$, a stable quasi-interpolant can be constructed via the formula
\begin{equation}
    z_{f,X} = \sum_{j=1}^{N} f(x_j)u_j.
    \label{eq_pol_rep_quasi_interpolant}
\end{equation}
\noindent It is worth noting that the second condition in Definition \ref{def_local_app_procc} ensures the stability of the quasi-interpolation process, since it states that the associated Lebesgue constant remains bounded. Indeed,
\begin{equation*}
	\begin{split}
		 \left|\sum_{j=1}^N f(x_j)u_j(x)-\sum_{j=1}^N \widetilde{f}(x_j)u_j(x)\right| & \leq \sum_{j=1}^N |f(x_j)-\widetilde{f}(x_j)||u_j(x)| \leq \\
		 & \leq \max_{1 \leq j \leq N} |f(x_j)-\widetilde{f}(x_j)| \sum_{j=1}^N |u_j(x)| \leq \\
		 & \leq \max_{1 \leq j \leq N} |f(x_j)-\widetilde{f}(x_j)| C_1.
	\end{split}
\end{equation*}
Instead, the first and third conditions are directly related to local polynomial reproduction. In this setting, the functions $\{u_j\}_{j=1,\dots,N}$ can be chosen quite arbitrarily, as their smoothness does not affect the stability or convergence of the quasi-interpolation process. Although the fast-decaying polynomial reproduction framework is quite general, it is beneficial to impose additional regularity assumptions on both the domain and the node distribution in order to derive concrete approximation schemes. To this end, we introduce two well-known definitions: the interior cone condition for the approximation domain and the notion of quasi-uniform node distribution.
\begin{definition}
A bounded Lipschitz domain $\Omega \subseteq \mathbb{R}^d$  satisfies an {\it interior cone condition}, if there exists an angle $\vartheta \in (0,\pi/2)$ and a radius $r>0$ such that for every $x \in \Omega$ a unit vector $\xi(x)$ exists so that the cone
	\begin{equation*}
		C(x,\xi(x),\vartheta,r) \coloneqq \{ x+ \lambda y : y \in \mathbb{R}^d,\, \|y\|_2 = 1, \, \langle y, \xi(x) \rangle \geq \cos (\vartheta), \, \lambda \in [0,r] \}
	\end{equation*}
 is entirely contained in $\Omega$. 
\end{definition}
\begin{definition}\label{def_quasi_uniform_data_set}
A nodes distribution $X=\{x_1,\dots,x_N\}$ is quasi-uniform with respect to a constant $c_{qu}>0$, if
\begin{equation*}
        q_{X} \leq h_{X,\Omega} \leq c_{qu} q_{X}\,.
\end{equation*}
\end{definition}
\noindent The requirement $q_{X} \leq h_{X,\Omega}$ instead of $c q_{X} \leq h_{X,\Omega}$ with $c < 1$ is not restrictive if $\Omega$ satisfies an interior cone condition with radius $r>0$ and angle $\vartheta>0$. In fact, if $q_{X} \leq r$ then for each $x_i \in X$ we can find $y \in \Omega$ with $\|y-x_i\|_2 = q_{X}$. Indeed, if $\|z\|_2=1$ and $\langle z, \xi(x_i) \rangle \geq \cos(\vartheta)$ then
\begin{equation*}
    \|\underbrace{x_i+\lambda z}_{y} -x_i\|_2 = \lambda \in [0,r]\,.
\end{equation*}
For any other $x_j \in \Omega$, we have
\begin{equation*}
\|y-x_j\|_2 \geq \|x_j-x_i\|_2-\|y-x_i\|_2 \geq 2q_{X}-q_{X}=q_{X}.
\end{equation*}
This implies $h_{X,\Omega} \geq q_{X}$.
\\ \\
The definition \ref{def_quasi_uniform_data_set} is particularly useful when considering a sequence of data sets sharing the same constant $c_{qu}$ while the fill distance $h_{X,\Omega}$ tends to zero.
\\
In the context of the RL-RBF (cf. Section \ref{sec_stability_convergence}) approximation with quasi-uniform node distributions, it is essential to control the number of nodes of $X$ lying near each point $x \in \Omega$ (we are analyzing the so-called \textit{stationary approach}, cf. \cite[Theorem 11.17]{wendland_2004}). This is achieved by ensuring that the scaling parameter $\delta$ of the approximation scheme, which corresponds to the reciprocal of the \textit{shape parameter} employed in RBF interpolation, remains bounded in terms of the fill distance (cf. \cite{demarchi2020})
 \begin{equation} 
    \gamma c_{\gamma} h_{X,\Omega} \leq \delta \leq c_{\gamma} h_{X,\Omega},
    \label{eq_delta_h_X_relation}
\end{equation}
with $\gamma \in ]0,1[$ and $c_{\gamma}>1$. 
If $X$ is a quasi-uniform data distribution, inequality \eqref{eq_delta_h_X_relation} becomes
\begin{equation}
    \gamma c_{\gamma} q_{X} \leq \delta \leq c_\gamma c_{qu} q_X.
    \label{eq_delta_q_X_relation}
\end{equation}
Finally, the existence of methods that provide local polynomial reproduction is guaranteed by the following result \cite[Theorem 3.14]{wendland_2004}.
\begin{theorem}
	Suppose that $\Omega \subseteq \mathbb{R}^d$ is compact and satisfies an interior cone condition with angle $\vartheta \in ]0,\pi/2[$ and radius $r>0$. Fix $m \in \mathbb{N}$. Then, there exist constants $h_0, C_1, C_2$ depending on $m, \vartheta,r$ such that for every $X=\{x_1, \dots, x_N\} \subseteq \Omega$ with $h_{X,\Omega} \leq h_0$ and every $x \in \Omega$ we can find real numbers $\{\widetilde{u}_j(x)\}_{j=1,\dots,N}$ with
	\begin{itemize}
		\item $\displaystyle \sum_{j=1}^{N} p(x_j)\widetilde{u}_j(x) = p(x)$ for each $p \in \pi_{m}(\mathbb{R}^d)$,
		\item $\displaystyle \sum_{j=1}^{N} |\widetilde{u}_j(x)| \leq C_1$,
		\item $\widetilde{u}_j(x)=0$ if $\|x-x_j\|_2 > C_2 h_{X,\Omega}$.
	\end{itemize}
 \label{thm_constant_lpr}
\end{theorem}
\noindent In particular, the proof of the preceding theorem yields the following explicit values for the constants involved, provided that $\vartheta \leq \pi/5$:
\begin{equation*}
	C_1 = 2, \qquad C_2 = \frac{16(1+\sin(\vartheta))^2m^2}{3\sin(\vartheta)^2}, \qquad h_0 = \frac{r}{C_2}.
\end{equation*}

\section{Fast-decaying polynomial reproduction}
\label{sec_stability_convergence}
Generalizing Definition \ref{def_local_app_procc}, we aim to relax the assumption of compact support while preserving the same rate of convergence (cf. \cite[Theorem 3.2]{wendland_2004}).
\begin{definition}
Let  $\varphi : \mathbb{R}_{\geq 0} \rightarrow \mathbb{R}$ be a decreasing function such that $\lim_{n\to +\infty} \frac{\varphi(n+1)}{\varphi(n)}$ exists and it is strictly smaller than $1$. A numerical scheme that produces for every set $X=\{ x_1, \dots, x_N\} \subseteq \Omega$ a family of functions $u_j = u_j^X : \Omega \rightarrow \mathbb{R}, \, 1 \leq j \leq N$, provides fast-decaying polynomial reproduction of degree $\ell$ on $\Omega$ with respect to $\varphi$ if there exist constants $C,h_0 \in \mathbb{R}_{>0}$ such that
	\begin{enumerate}
		\item $\sum_{j=1}^{N} p(x_j)u_j = p$ for each $p \in \pi_{\ell}(\mathbb{R}^d)$,
		\item $|u_j(x)| \leq C\varphi\left(\frac{\|x-x_j\|_2}{q_X}\right)$ for all $x \in \Omega$ and $j=1,\dots,N$
	\end{enumerate}
are satisfied for all $X$ with $h_{X,\Omega} \leq h_0$.
\label{def_fast_decay_pol_rep}
\end{definition}
\noindent In order to demonstrate the utility of the  definition \ref{def_fast_decay_pol_rep}, we apply it to the RL-RBF method. Let us start by fixing a positive definite function $\Phi : \mathbb{R}^d \rightarrow \mathbb{R}$ with support in the closure of the unit ball $B(0,1)$, such that its native space coincides with a Sobolev space. For $\delta>0$ (cf. eq. \eqref{eq_delta_h_X_relation}), we define the scaled function
\begin{equation*}
    \Phi_\delta(\cdot) = \Phi \left( \frac{\cdot}{\delta}\right).
\end{equation*}
Given interpolation nodes $X$ and a function $f \in \mathcal{C}(\Omega)$, one can construct the classical RBF interpolant:
\begin{equation} \label{RBF_RLRBF_interp_problem}
    s_{f,X,\Phi_\delta}(x) = \sum_{j=1}^N \alpha_j \Phi_\delta(x-x_j), \qquad x \in \Omega.
\end{equation}
Consequently, the RL-RBF interpolant takes the form
\begin{equation} \label{RL_RBF_approximant_definition_001}
    z_{f,X}=\frac{ s_{f,X,\Phi_\delta}(x)}{ s_{1,X,\Phi_\delta}(x)}.
\end{equation}
If $\{\chi_1, \dots, \chi_N\}$ are the cardinal functions associated with the interpolation problem in equation \eqref{RBF_RLRBF_interp_problem}, then the RL-RBF interpolant becomes
\begin{equation*}
    z_{f,X}(x)= \frac{\sum_{j=1}^N f(x_j) \chi_j(x)}{\sum_{k=1}^N \chi_k(x)}= \sum_{j=1}^N f(x_j)u_j(x),
\end{equation*}
with the new basis functions
\begin{equation*}
    u_j(x) = \frac{\chi_j(x)}{\sum_{k=1}^N \chi_k(x)}.
\end{equation*}
In this context, for the sake of brevity we do not discuss the RL-RBF method further; however, by means of \cite[Theorem 2.3, Lemma 3.1]{demarchi2020} we can show that the method just described fits naturally within the framework of Definition \ref{def_fast_decay_pol_rep} as it exactly reproduces constant functions,
\begin{equation*}
    |u_j(x)| \leq C \exp\left({-\frac{\nu \| x-x_j\|_2}{q_X}}\right), \qquad x \in \mathbb{R}^d.
\end{equation*}
In what follows, we address the problem of constructing a stable approximation based on the fast-decaying polynomial reproduction scheme, as given in eq. \eqref{eq_pol_rep_quasi_interpolant}, for any function $f: \Omega^* \rightarrow \mathbb{R}$ which is ${\cal C}^{m+1}(\Omega^*)$, where $\Omega^*$ denotes the closure of the convex hull of $\Omega$ (see Theorem \ref{thm_convergence_polynomial_reproduction_exponential_decay} below).

\subsection{Stability}
The following result, which is a straightforward generalization of \cite[Corollary 2.5]{demarchi2020}, establishes the stability of the fast-decaying polynomial reproduction framework and will be instrumental in proving its convergence.
\begin{theorem}
    Let $X=\{x_1,\dots,x_N\} \subseteq \mathbb{R}^d$ be an arbitrary data set, and let $\{u_1,\dots,u_N\}$ be given by a fast-decaying polynomial reproduction method with respect to $\varphi$. Then, for every $\ell \in \mathbb{N}$ there exists a constant $K=K(\ell,C,\varphi,d)$ such that 
    \begin{equation}
        \sum_{j=1}^{N} \|x-x_j\|_2^{\ell} |u_j(x)| \leq K q_{X}^{\ell} \quad \textrm{ for } x \in \Omega,
         \label{eq_stability_quasi_interpoaltion_fast_decay}
    \end{equation}
    where
    \begin{equation}
        K=3^d C \sum_{n=0}^{+\infty} (n+1)^{d+\ell-1}\varphi(n).
        \label{def_constant_K_stability}
    \end{equation}
    \label{thm_stability_quasi_interpoaltion_fast_decay}
\end{theorem}
\begin{proof}
Fix $x \in \Omega$. We define a sequence of sets $\{E_n\}_{n \in \mathbb{N}}$ that covers $\mathbb{R}^d$ as
\begin{equation*}
    E_n = \{ y \in \mathbb{R}^d: n q_X \leq \|y-x\|_2 \leq (n+1)q_X\}.
\end{equation*}
We note that if $x_j \in E_n$ then
\begin{equation*}
    B(x_j,q_X) \subseteq \{ y \in \mathbb{R}^d: (n-1)q_X \leq \|y-x\|_2 \leq (n+2)q_X\}=E_{n-1}\cup E_n \cup E_{n+1}\,.
\end{equation*}
In fact, if $y \in B(x_j,q_X)$ then $nq_X -q_X \leq \|x_j-x\|_2-\|x_j-y\|_2 \leq\|y-x\|_2 \leq \|y-x_j\|_2+\|x_j-x\| \leq q_X +(n+1)q_X$.
Since $\{B(x_j,q_X)\}_{j=1,\dots,N}$ are pairwise disjoint,
\begin{equation*}
    \bigcup_{j : x_j \in E_n} B(x_j,q_X) \subseteq \overline{B(x,(n+2)q_X)} \setminus B(x,(n-1)q_X),
\end{equation*}
a volume comparison gives for $n \geq 1$
\begin{equation*}
    \# \{x_j : x_j \in E_n\} \leq (n+2)^d-(n-1)^d \leq 3^d n^{d-1} \leq 3^d(n+1)^{d-1},
\end{equation*}
that holds also for $n=0$ because $2^d \leq 3^d$ ($\#$ being the cardinality of the set).
\\
\\
The foregoing identity is established by induction on the dimension. For $d=1$ we have $n+2-n+1=3$. 
\\ Let us suppose that it is true for $d$, we will prove the inequality for $d+1$.
\begin{equation*}
    \begin{split}
        &(n+2)^{d+1}-(n-1)^{d+1} = (n+2)(n+2)^d-(n-1)(n-1)^d = \\
        = & n((n+2)^d-(n-1)^d)+2(n+2)^d+(n-1)^d = \\
        = & n((n+2)^d-(n-1)^d) +2((n+2)^d-(n-1)^d)+3(n-1)^d \leq\\
        \leq & 3^dn^d+2 \cdot3^dn^{d-1}+3(n-1)^d \leq 3^{d+1}n^d,
    \end{split}
\end{equation*}
where the last inequality holds because by dividing by $3^dn^d$ we have
\begin{equation*}
    1+\frac{2}{n}+ \frac{3}{3^d}\left(\frac{n-1}{n}\right)^d \leq 3,
\end{equation*}
that for $n=1$ becomes $3 \leq 3$, instead for $n\geq 2$ is 
\begin{equation*}
     \frac{2}{n}+ \frac{3}{3^d}\left(\frac{n-1}{n}\right)^d \leq \frac{2}{n}+\left(\frac{n-1}{n}\right)^d \leq \frac{2}{n}+1 \leq 1+1 = 2.
\end{equation*}
Observing that if $x_j \in E_n$ then $n \leq \frac{\|x_j-x\|_2}{q_X}$ and if we split the sum over the sets $\{E_n\}_{n \in \mathbb{N}}$ we obtain 
\begin{equation*}
\begin{split}
    \sum_{j=1}^{N} \|x-x_j\|_2^\ell|u_j(x)| & \leq \sum_{n=0}^{+\infty} \sum_{x_j \in E_n} \|x-x_j\|_2^\ell|u_j(x)| \leq \sum_{n=0}^{+\infty}\sum_{x_j \in E_n} \|x-x_j\|_2^\ell C\varphi\left( \frac{\|x-x_j\|_2}{q_X} \right) \leq \\
    & \leq \sum_{n=0}^{+\infty}\sum_{x_j \in E_n} \|x-x_j\|_2^\ell C\varphi(n) \leq \sum_{n=0}^{+\infty} \sum_{x_j \in E_n} (n+1)^\ell q_X^\ell C\varphi(n) \leq \\
    & \leq \sum_{n=0}^{+\infty} 3^d (n+1)^{d+\ell-1}q_X^\ell C \varphi(n) = q_{X}^{\ell} \underbrace{3^d C \sum_{n=0}^{+\infty} (n+1)^{d+\ell-1}\varphi(n)}_{K}.
\end{split}
\end{equation*}
Convergence of the series follows from the ratio test:
\begin{equation*}
    \frac{(n+2)^{d+\ell-1} \varphi(n+1)}{(n+1)^{d+\ell-1} \varphi(n)} = \left( \frac{n+2}{n+1}\right)^{d+\ell-1} \frac{\varphi(n+1)}{\varphi(n)} \xrightarrow[]{n \rightarrow +\infty}  \lim_{n\to +\infty} \frac{\varphi(n+1)}{\varphi(n)} < 1.
\end{equation*}
\end{proof}
\noindent We remark that quasi-uniformity of the set $X$ is not required for the present proof, although it will play a crucial role in the convergence analysis. By Theorem \ref{thm_stability_quasi_interpoaltion_fast_decay}, the Lebesgue function associated with the quasi-interpolation scheme \eqref{eq_pol_rep_quasi_interpolant} satisfies
\begin{equation*}
     \sum_{j=1}^N |u_j(x)| \leq  3^d C \sum_{n=0}^{+\infty} (n+1)^{d-1}\varphi(n)
\end{equation*}
uniformly for all $x\in \Omega$. This bound establishes the stability of the scheme.
\\
\\
The decay conditions in Definition \ref{def_fast_decay_pol_rep} are formulated as geometric (exponential) decay conditions, which facilitates the proof of Theorem \ref{thm_stability_quasi_interpoaltion_fast_decay}. No specific decay rate is required to establish the result, since the only essential condition is the convergence of the series appearing in equation \eqref{def_constant_K_stability}. In Section \ref{subsection_algebraic_decay_examples}, we shall demonstrate how the fast-decaying polynomial reproduction framework can be applied to degenerate cases involving algebraic decay.

\subsection{Convergence}
The framework under consideration requires the exact reproduction of polynomial spaces. Consequently, if all polynomials of degree at most $m$ are reproduced exactly, the approximation is expected to converge with order $\mathcal{O}(h_{X,\Omega}^{m+1})$.
\begin{theorem}
Let $X$ and $\{u_j\}_{j=1}^N$ as in the Theorem \ref{thm_stability_quasi_interpoaltion_fast_decay}, and assume in addition that $X$ is quasi-uniform.
    %Suppose that $X=\{x_1,\dots,x_N\} \subseteq \Omega$ is an arbitrary data set and $\{u_1,\dots,u_N\}$ is given by a fast decaying polynomial reproduction process with respect to $\varphi$ of order $m \in \mathbb{N}$. 
    Let $\Omega^{*}$ denote the closure of the convex hull of $\Omega$, and assume that $ \Omega $ satisfies an interior cone condition. If $f \in \mathcal{C}^{m+1}(\Omega^{*})$, then there exists a constant $K=K(C,\varphi,d,m)$ such that 
    \begin{equation*}
        \|f-z_{f,X}\|_{L^{\infty}(\Omega)} \leq K h_{X,\Omega}^{m+1}\|f\|_{\mathcal{C}^{m+1}(\Omega^*)},
    \end{equation*}
    for $h_{X,\Omega} \leq h_0$.
    \label{thm_convergence_polynomial_reproduction_exponential_decay}
\end{theorem}
\begin{proof}
Let $p \in \pi_{m}(\mathbb{R}^d)$ be an arbitrary polynomial of degree $\le m$. By means of Definition \ref{def_fast_decay_pol_rep} 
\begin{equation*}
	\begin{split}
		|f(x)-z_{f,X}(x)| & \leq |f(x)-p(x)|+\left| p(x)-\sum_{j=1}^N f(x_j)u_j(x)\right| \leq \\
		& \leq |f(x)-p(x)| + \sum_{j=1}^N |p(x_j)-f(x_j)||u_j(x)|\,.
  \end{split}
\end{equation*}
Choose $p$ to be the Taylor polynomial of $f$ around $x$ of order $m$, then for $y \in \Omega$ there exists $\xi \in \Omega^{*}$ such that
\begin{equation*}
	f(y)-\sum_{|\alpha| \leq m} \frac{D^{\alpha}f(x)}{\alpha!}(y-x)^{\alpha} = \sum_{|\alpha|= m+1} \frac{D^{\alpha}f(\xi)}{\alpha!}(y-x)^{\alpha}.
\end{equation*}
Remarking that $|(y-x)^\alpha| = \prod_{i=1}^{d}|y_i-x_i|^{\alpha_i} \leq \prod_{i=1}^{d}\|y-x\|^{\alpha_i} = \|y-x\|^{|\alpha|}$, for any $j=1,\dots,N$
\begin{equation*}
    |p(x_j)-f(x_j)| \leq \sum_{|\alpha|=m+1}\frac{\|D^{\alpha}(f)\|_{L^{\infty}(\Omega^*)}}{\alpha!}\|x_j-x\|_2^{m+1} \leq \left( \sum_{|\alpha|=m+1}\frac{1}{\alpha!} \right)|f|_{\mathcal{C}^{m+1}(\Omega^*)}\|x_j-x\|_2^{m+1}\,.
\end{equation*}
Hence,
\begin{equation*}
\begin{split}
    |f(x)-z_{f,X}(x)| & \leq \left( \sum_{|\alpha|=m+1}\frac{1}{\alpha!} \right)|f|_{\mathcal{C}^{m+1}(\Omega^*)} \sum_{j=1}^{N} \|x_j-x\|_2^{m+1} |u_{j}(x)| \\ 
    & \leq \left( \sum_{|\alpha|=m+1}\frac{1}{\alpha!} \right) K(m+1,C,\varphi,d) |f|_{\mathcal{C}^{m+1}(\Omega^*)} h_{X,\Omega}^{m+1},
\end{split}
\end{equation*}
where $K(m+1,C,\varphi,d)$ comes from the application of Theorem \ref{thm_stability_quasi_interpoaltion_fast_decay} with $q_X \leq h_{X,\Omega}$ and the fact that $\Omega$ satisfies an interior cone condition.
\end{proof}

\section{An example of fast-decaying polynomial reproduction}
\label{section_ex_c_infinity_fast}

The fast-decaying polynomial reproduction framework can be naturally incorporated into the moving least squares method \cite{lancaster_1981,mclain_1974,mclain_1976,shepard_1968}. In particular, by employing Gaussian weights, it is possible to construct a $ C^\infty(\Omega)$ approximant. For completeness, we recall that the moving least squares approximant $z_{f,X}(x)$ at any point $x \in \Omega$ is given by $z_{f,X}(x) = p^*(x)$, where $p^{*}$ is the unique solution of the minimization problem
\begin{equation}
	\min \left\{ \sum_{i=1}^{N}(f(x_i)-p(x_i))^2e^{-\nu \left(\frac{\|x-x_i\|_2}{\delta}\right)^2} : p \in \pi_{m}(\mathbb{R}^d)\right\},
 \label{eq_moving_least_squares_fast_decay}
\end{equation}
with $\nu \in \mathbb{R}_{>0}$ a fixed parameter and $\delta$ (cf. eq. \eqref{eq_delta_q_X_relation}) depending on the data set $X$.
\\
\\
Since the weight function 
\begin{equation*}
    w(x,y) = e^{-\nu \left(\frac{\|x-y\|_2}{\delta}\right)^2}>0
\end{equation*} is strictly positive, we can generalize the results of \cite[Section 4.1]{wendland_2004}. Indeed, if the set $X=\{x_1,\dots,x_N \}$ is $\pi_{m}(\mathbb{R}^d)$-unisolvent, then the unique solution of \eqref{eq_moving_least_squares_fast_decay} is
\begin{equation*}
	z_{f,X}(x) = \sum_{i=1}^{N} f(x_i)a_i^{*}(x),
\end{equation*}
where the coefficient vector $\boldsymbol{a}^*(x) = (a_1^*(x), \dots, a_N^*(x))^T$ is the unique minimizer of the constrained quadratic optimization problem
\begin{equation*}
\min_{\boldsymbol{a}(x) \in \mathbb{R}^N} \boldsymbol{a}^T(x)D(x) \boldsymbol{a}(x),
\end{equation*}
with $D(x)$ a diagonal matrix with elements $d_i(x)=e^{-\nu \left(\frac{\|x-x_i\|_2}{\delta}\right)^2}$, that is  
\begin{equation}
	\min_{\boldsymbol{a}(x) \in \mathbb{R}^N} \sum_{i=1}^N \frac{a_{i}(x)^2}{e^{-\nu \left(\frac{\|x-x_i\|_2}{\delta}\right)^2}}
	\label{eq_def_minimization_mls_2_exp_fast_decay}
\end{equation}
under the constraints 
\begin{equation}
	\sum_{i=1}^N p(x_i)a_{i}(x) = p(x), \qquad p \in \pi_{m}(\mathbb{R}^d).
	\label{eq_def_minimization_mls_2_const_exp_fast_decay}
\end{equation}
\noindent It is also possible to derive an explicit representation for the basis functions $\{a_{j}^{*}(x)\}_{j=1,\dots,N}$ as
\begin{equation*}
    a_{j}^{*}(x) = e^{-\nu \left(\frac{\|x-x_j\|_2}{\delta}\right)^2}\sum_{k=1}^{Q}\lambda_k(x)p_k(x_j),
\end{equation*}
where $\{\lambda_k(x)\}_{k=1,\dots,Q}$ are the unique solutions of
\begin{equation}
    \sum_{k=1}^{Q} \lambda_k(x) \sum_{j=1}^{N} e^{-\nu \left(\frac{\|x-x_j\|_2}{\delta}\right)^2}p_k(x_j)p_{\ell}(x_j)=p_{\ell}(x), \qquad 0 \leq \ell \leq Q.
    \label{eq_del_lambda_k_x_exp_fast_decay}
\end{equation}
 Hence, being $\{p_1, \dots, p_Q\}$ a basis for $\pi_{m}(\mathbb{R}^d)$, it follows that the approximant $z_{f,X} \in \mathcal{C}^{\infty}(\Omega)$. 
 \\ \\
For the case $m=0$, which corresponds to the reproduction of constants, we obtain
\begin{equation*}
    a_{j}^{*}(x) = \frac{e^{-\nu \left(\frac{\|x-x_j\|_2}{\delta}\right)^2}}{\sum_{i=1}^{N} e^{-\nu \left(\frac{\|x-x_i\|_2}{\delta}\right)^2}}\,, \qquad  j =1, \dots, N,
\end{equation*}
and the minimization problem \eqref{eq_moving_least_squares_fast_decay} has the solution 
\begin{equation}
    z_{f,X}(x) = \sum_{j=1}^{N} f(x_j) \underbrace{\frac{e^{-\nu \left(\frac{\|x-x_j\|_2}{\delta}\right)^2}}{\sum_{i=1}^N e^{-\nu \left(\frac{\|x-x_i\|_2}{\delta}\right)^2}}}_{a_{j}^{*}(x)},
    \label{eq_shepard_method_general}
\end{equation}
which constitutes a particular instance of the Shepard approximation method \cite{shepard_1968}.
\\
A similar result for compactly supported basis functions in quasi-uniform settings is established in \cite[Theorem 4.7]{wendland_2004}. In the present work, we extend this result to the more general setting of fast-decaying basis functions.
\begin{theorem}
    Assume that $\Omega \subseteq \mathbb{R}^d$ is compact and satisfies an interior cone condition with angle $\vartheta \in ]0,\pi/2[$ and radius $r>0$. Let $m \in \mathbb{N}$ be fixed. Denote $h_0, C_1$ and $C_2$ the constants appearing in Theorem \ref{thm_constant_lpr}. Suppose that $X=\{x_1, \dots, x_N\} \subseteq \Omega$ is a quasi-uniform data sets with constant $c_{qu}>0$ satisfying $h_{X,\Omega} \leq h_0$. Let $\delta$ be as in equation \eqref{eq_delta_h_X_relation}. Then the basis functions $\{a_{j}^{*}(x)\}_{j=1,\dots,N}$ provide fast-decaying polynomial reproduction, where the constants $C$, $h_0$ and the decay function $\varphi$ can be determined explicitly.
    \label{thm_pr_fast_decay_exp}
\end{theorem}

\begin{proof}
By Definition \ref{def_fast_decay_pol_rep} we must prove Properties 1 and 2.
The first property of Definition \ref{def_fast_decay_pol_rep} is a consequence of equation \eqref{eq_moving_least_squares_fast_decay} that defines the moving least squares method.
\\
\\
To prove the second assertion, we need to bound the quantity
\begin{equation*}
\begin{split}
     \sum_{i=1}^{N} \frac{(a_{i}^{*}(x))^2}{e^{-\nu \left(\frac{\|x-x_i\|_2}{\delta}\right)^2}}.
\end{split}
\end{equation*}
There exists $\{\widetilde{u}_j(x)\}_{j=1,\dots,N}$ providing local polynomial reproduction (cf. Theorem \ref{thm_constant_lpr}) such that $\widetilde{u}_j$ is supported in $\overline{B(x_j, C_2 h_{X,\Omega})}$ for $j=1,\dots,N$.
Being basis function $\{a_{j}^{*}(x)\}_{j=1,\dots,N}$ obtained by solving the minimization problem \eqref{eq_def_minimization_mls_2_const_exp_fast_decay}, we get
\begin{equation}
    \begin{split}
        \sum_{i=1}^N \frac{(a_{i}^{*}(x))^2}{e^{-\nu \left(\frac{\|x-x_i\|_2}{\delta}\right)^2}} & \leq \sum_{i=1}^{N} \frac{(\widetilde{u}_i(x))^2}{e^{-\nu \left(\frac{\|x-x_i\|_2}{\delta}\right)^2}} = \sum_{i \in \widetilde{I}(x)} \frac{(\widetilde{u}_i(x))^2}{e^{-\nu \left(\frac{\|x-x_i\|_2}{\delta}\right)^2}} \leq \\
        & \leq \frac{1}{e^{-\nu \left(\frac{C_2 h_{X,\Omega}}{\delta}\right)^2}} \sum_{i \in \widetilde{I}(x)} (\widetilde{u}_i(x))^2 \leq \\
        & \leq \frac{1}{e^{-\nu \left(\frac{C_2 h_{X,\Omega}}{\delta}\right)^2}} \left( \sum_{i \in \widetilde{I}(x)} |\widetilde{u}_i(x)| \right)^2 \leq \\
        & \leq \frac{C_1^2}{e^{-\nu \left(\frac{C_2 h_{X,\Omega}}{\delta}\right)^2}} \leq \frac{C_1^2}{e^{-\nu \left(\frac{C_2}{\gamma c_\gamma}\right)^2}},
    \end{split}
    \label{eq_computation_proof_mls}
\end{equation}
where $\widetilde{I}(x)=\left\{j \in \{1,\dots,N\} : x_j \in \overline{B(x,C_2 h_{X,\Omega})}\right\}$. With the last computation we obtained for $i=1,\dots,N$
\begin{equation*}
    \frac{(a_{i}^{*}(x))^2}{e^{-\nu \left(\frac{\|x-x_i\|_2}{\delta}\right)^2}} \leq \frac{C_1^2}{e^{-\nu \left(\frac{C_2}{\gamma c_\gamma}\right)^2}} \Rightarrow (a_{i}^{*}(x))^2 \leq \frac{C_1^2}{e^{-\nu \left(\frac{C_2}{\gamma c_\gamma}\right)^2}} e^{-\nu \left(\frac{\|x-x_i\|_2}{\delta}\right)^2}
\end{equation*}
that gives us
\begin{equation*}
    |a_{i}^{*}(x)| \leq \sqrt{\frac{C_1^2}{e^{-\nu \left(\frac{C_2}{\gamma c_\gamma}\right)^2}} } e^{-\frac{\nu}{2} \left(\frac{\|x-x_i\|_2}{\delta}\right)^2} \leq \sqrt{\frac{C_1^2}{e^{-\nu \left(\frac{C_2}{\gamma c_\gamma}\right)^2}}} e^{-\frac{\nu}{2} \left(\frac{\|x-x_i\|_2}{c_\gamma c_{qu} q_X}\right)^2}.
\end{equation*}
We now observe that $e^{-x^2} \leq e\, e^{-x}$ for $x \in \mathbb{R}$. This results from the fact that
\begin{equation*}
    -x^2 \leq -x +1 \Longleftrightarrow x \leq x^2 +1
\end{equation*}
that holds for $x \in [-1,1]$ and $x \in \mathbb{R} \setminus [-1,1]$. We then obtain
\begin{equation}
    |a_{i}^{*}(x)| \leq \underbrace{e \sqrt{\frac{C_1^2}{e^{-\nu \left(\frac{C_2}{\gamma c_\gamma}\right)^2}}}}_{C} e^{-\sqrt{\frac{\nu}{2}} \frac{\|x-x_i\|_2}{c_\gamma c_{qu} q_X}} = C\varphi \left( \frac{\|x-x_i\|_2}{q_X}\right)
    \label{eq_computation_proof_mls_1}
\end{equation} 
for $i=1,\dots,N$ with
\begin{equation*}
    \varphi(x) = e^{-\sqrt{\frac{\nu}{2}} \frac{1}{c_\gamma c_{qu}}x}.
\end{equation*}
\end{proof}
\noindent We observe that in \cite[Theorem 4.7]{wendland_2004} the set $X$ is required to be locally unisolvent, which may lead to oversampling. In contrast, Theorem \ref{thm_pr_fast_decay_exp} above only requires $X$ to be (globally) unisolvent, since the proof does not rely on local arguments.
\\
\\
The same construction used to establish Theorem \ref{thm_pr_fast_decay_exp} can also be carried out with the weight function 
\begin{equation*}
    w(x,y) = e^{-\nu \frac{\|x-y\|_2}{\delta}}>0,
\end{equation*}
yielding the explicit constants
\begin{equation}
    C = \sqrt{\frac{C_1^2}{e^{-\nu \left(\frac{C_2}{\gamma c_\gamma}\right)}}} \quad \textrm{ and } \quad \varphi(x) = e^{-\frac{\nu}{2} \frac{1}{c_\gamma c_{qu}}x}.
    \label{eq_stability_fast_decay_choose_delta}
\end{equation}
More generally, Theorem \ref{thm_pr_fast_decay_exp} remains valid for weight functions of the form
\begin{equation*}
    w(x,y)=\varphi \left(\frac{\|x-y\|_2}{q_X}\right)>0,
\end{equation*}
where $\varphi$ satisfies the assumptions of Definition \ref{def_fast_decay_pol_rep}. In this case, the constants are given by
\begin{equation}
    C = \sqrt{\frac{C_1^2}{\varphi( C_2 c_{qu})} } \quad \textrm{ and } \quad \widetilde{\varphi}(x) = \sqrt{\varphi(x)}.
    \label{eq_general_phi_weight_mls}
\end{equation}

For the case $m=0$, the minimization problem admits the explicit solution
\begin{equation*}
    z_{f,X}(x) = \sum_{j=1}^{N} f(x_j) \underbrace{\frac{\varphi \left(\frac{\|x-x_j\|_2}{q_X}\right)}{\sum_{i=1}^N \varphi \left(\frac{\|x-x_i\|_2}{q_X}\right)}}_{a_{j}^{*}(x)}.
\end{equation*}
The relationship between local polynomial reproduction and moving least squares approximation was investigated in \cite{wendland_2001} using compactly supported basis functions. To emphasize the generality of our framework, we present additional kernels with exponential decay that are well suited to our purposes, for example Mat\'ern functions \cite{matern_2013}:
\begin{equation*}
    \varphi(x) = \frac{K_{\beta-d/2}(x)x^{\beta-d/2}}{2^{\beta-1}\Gamma(\beta)}, \qquad d < 2\beta,
\end{equation*}
where $K_p$ is the modified Bessel function of the second kind of order $p$.
\begin{table}[H]
\begin{center}
\begin{tabular}{Sc | Sc Sc Sc}
\hline
$\beta$ & $\frac{d+1}{2}$ & $\frac{d+3}{2}$ & $\frac{d+5}{2}$ \\ \hline
 $\varphi$ & $e^{-\|x\|}$ & $(1+\|x\|)e^{-\|x\|}$ & $(3+3\|x\|+3\|x\|^2)e^{-\|x\|}$ \\ \hline
   Smoothness & $\mathcal{C}(\mathbb{R}^d)$ & $\mathcal{C}^2(\mathbb{R}^d)$ & $\mathcal{C}^4(\mathbb{R}^d)$\\ \hline
\end{tabular}
\end{center}
\caption{Mat\'ern functions corresponding to different values of the parameter $\beta$ and relative smoothness.}
\end{table}
\noindent We also note that, in a quasi-uniform setting, the schemes of Definition \ref{def_local_app_procc} satisfy Definition \ref{def_fast_decay_pol_rep}, yielding
\begin{equation*}
    |u_{i}(x)| \leq \sum_{j=1}^N |u_{j}(x)| \leq C_1,\; \, i=1,\ldots, N\,.
\end{equation*}
Thus, if we choose $\widetilde{K}>0$ such that $C_1 \leq \widetilde{K} e^{-C_2}$, then for all $x\in \Omega \cap B(x_i, C_2 h_{X,\Omega})$ we obtain
\begin{equation*}
    |u_{i}(x)| \leq C_1 \leq \widetilde{K}e^{-C_2} \leq \widetilde{K} e^{-\frac{\|x-x_i\|_2}{h_{X,\Omega}}} \leq \widetilde{K}e^{-\frac{1}{c_{qu}} \frac{\|x-x_i\|_2}{q_X}}.
\end{equation*}
Regarding the computational cost of evaluating the approximant at a point $ x \in \Omega$, we observe that computing the coefficients $\{\lambda_1(x), \dots, \lambda_Q(x)\}$ requires solving a $ Q \times Q$ linear system at a cost of $ \mathcal{O}(Q^3)$. Assembling the matrix of the linear system has a cost of $\mathcal{O}(N Q^2)$, where $N$ denotes the number of data sites. Computing the basis functions $\{a_1^*(x), \dots, a_N^*(x)\}$ requires $\mathcal{O}(N Q)$ operations, since the evaluation of each basis function involves a number of multiplications and additions proportional to $Q$. In summary, the total cost of constructing the basis functions at a single point $x \in \Omega$ is
\begin{equation*}
\mathcal{O}(Q^3 + NQ^2 + NQ) = \mathcal{O}(N),
\end{equation*}
assuming that $Q$ (dimension of $\pi_{m}(\mathbb{R}^d)$) is fixed. Adding the $  \mathcal{O}(N)  $ cost required to evaluate the linear combination then gives an overall cost of $  \mathcal{O}(N)  $ per evaluation point. Consequently, evaluating the approximant at $ M $ distinct points incurs a total computational cost of $  \mathcal{O}(N M) $.

\section{Numerical experiments}
\label{sec_numerical_experiments}
\subsection{Basis Functions}

\begin{figure}[!h]
	\begin{subfigure}{0.30\textwidth}
		\includegraphics[width=\linewidth]{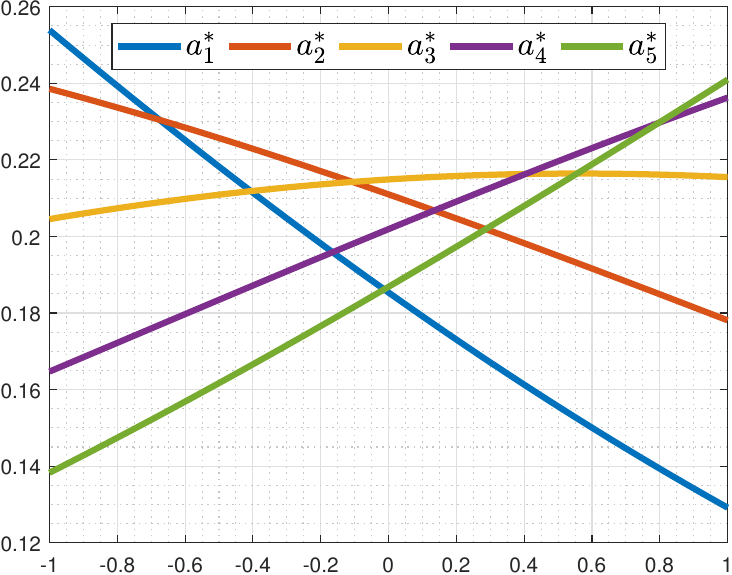}
		%		\label{fig:1}
  \caption{m=0}
  \end{subfigure}\hfil % <-- added
	\begin{subfigure}{0.30\textwidth}
		\includegraphics[width=\linewidth]{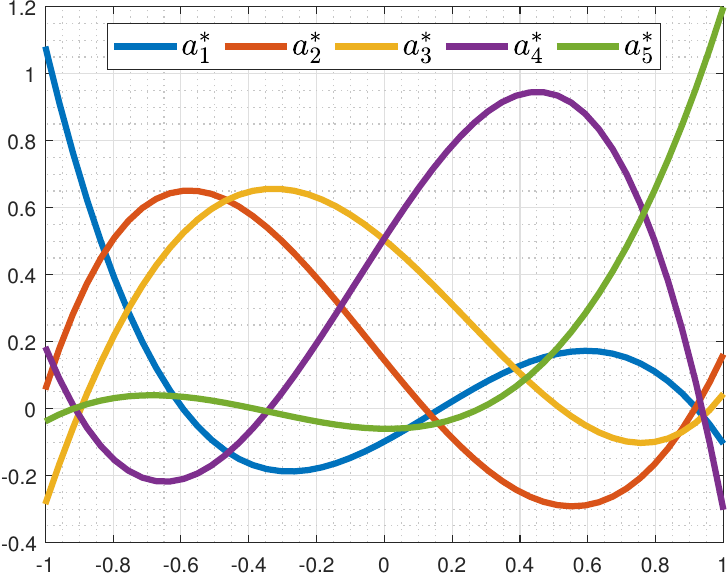}
		%		\label{fig:2}
  \caption{m=3}
  \end{subfigure}\hfil % <-- added
	\begin{subfigure}{0.30\textwidth}
		\includegraphics[width=\linewidth]{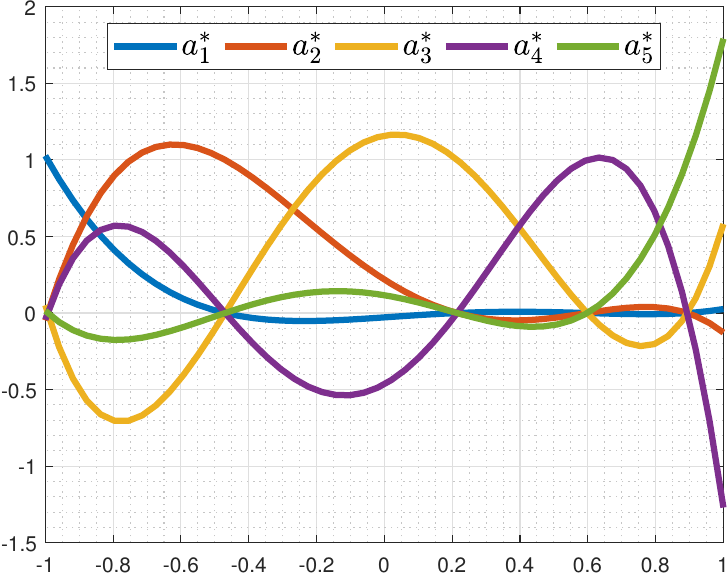}
		%		\label{fig:2}
  \caption{m=4}
  \end{subfigure}\hfil % <-- added
	
\caption{The weight function is given by $w(x)=e^{-x^2} \in \mathcal{C}^{\infty}(\mathbb{R})$, and the nodes consist of five uniformly perturbed equispaced points in the interval $[-1,1]$. From left to right, the basis functions reproduce polynomials of degrees $m=0, 3$, and $4$, respectively. In this experiment, the scaling parameter is set to $\delta = 5h_{X,\Omega}$.}
\label{fig_exp_2_basis_function_fast_decay}
\end{figure}

The purpose of this section is to present numerical experiments that illustrate the theoretical properties of fast-decaying polynomial reproduction methods, as introduced in Definition \ref{def_fast_decay_pol_rep}.
% and of the approximation scheme based on the $\|\cdot\|_1$ norm (see equation \eqref{approximant_linear_problem_optimization_form}). 
In contrast to previous works on compactly supported basis functions (e.g., \cite{fasshauer_2007}), we focus on non-compactly supported basis functions that enable the construction of smooth approximants.
\\
\\
We begin by examining the basis functions arising from Theorem \ref{thm_pr_fast_decay_exp}. Figure \ref{fig_exp_2_basis_function_fast_decay} illustrates the differentiability properties of these basis functions, in accordance with Corollary 4.5 of \cite{wendland_2004}. In particular, the basis functions are smooth regardless of the polynomial space they reproduce.

\begin{figure}[!h]
	\begin{subfigure}{0.30\textwidth}
		\includegraphics[width=\linewidth]{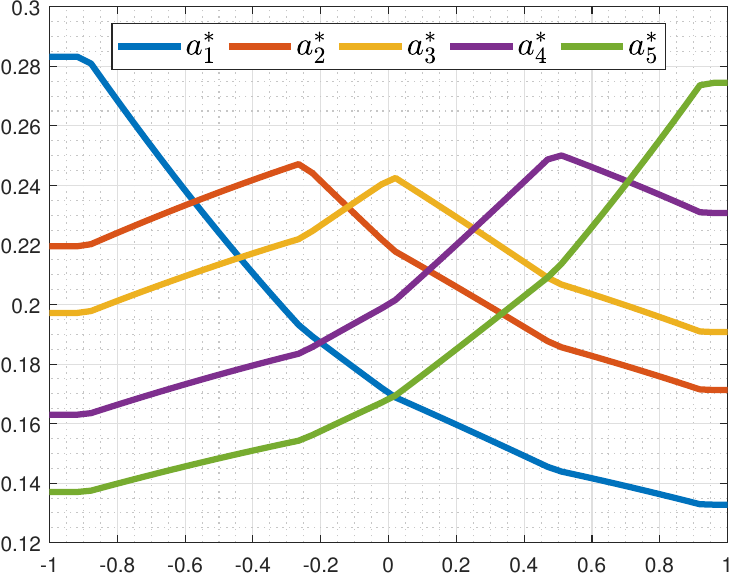}
		%		\label{fig:1}
  \caption{m=0}
  \end{subfigure}\hfil % <-- added
	\begin{subfigure}{0.30\textwidth}
		\includegraphics[width=\linewidth]{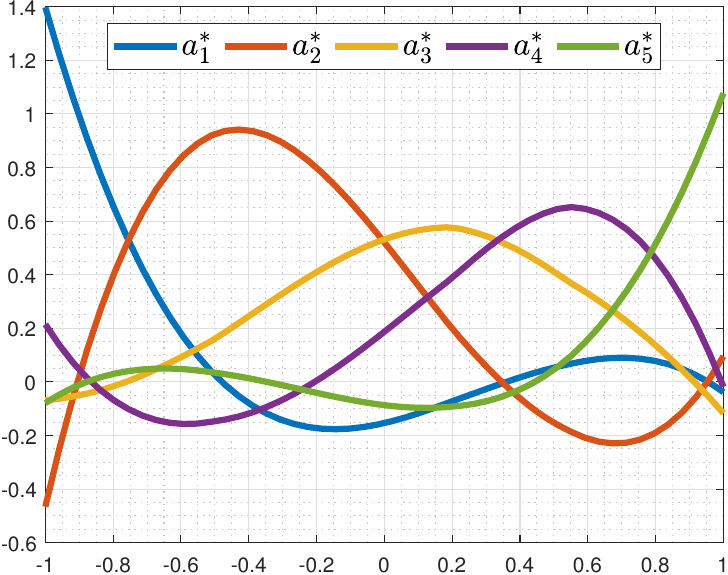}
		%		\label{fig:2}
  \caption{m=3}
  \end{subfigure}\hfil % <-- added
	\begin{subfigure}{0.30\textwidth}
		\includegraphics[width=\linewidth]{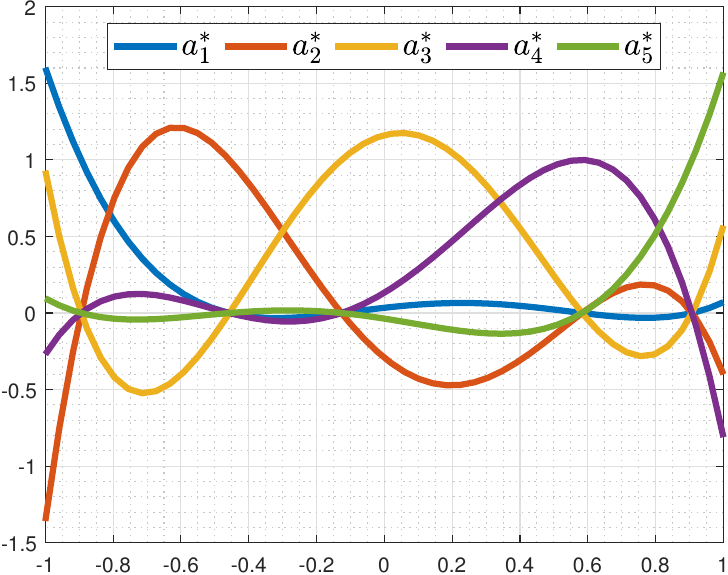}
		%		\label{fig:2}
  \caption{m=4}
  \end{subfigure}\hfil % <-- added
	
\caption{Basis functions arising from Theorem \ref{thm_pr_fast_decay_exp}. The weight function is given by $w(x)=e^{-x} \in \mathcal{C}(\mathbb{R})$, and the nodes consist of five uniformly perturbed equispaced points in the interval $[-1,1]$. From left to right, the basis functions reproduce polynomials of degrees $0, 3$ and $4$, respectively. In this experiment, the scaling parameter is set to $\delta = 5h_{X,\Omega}$.}
\label{fig_exp_basis_function_fast_decay}
\end{figure}

\vskip 0.1in
\noindent Figure \ref{fig_exp_basis_function_fast_decay} highlights an interesting feature regarding the smoothness of the basis functions. According to \cite[Corollary 4.5]{wendland_2004}, the basis functions $\{a^{*}_{1},\dots,a^{*}_{5}\}$ are a priori guaranteed to belong to $\mathcal{C}([-1,1])$, and this is confirmed by the numerical results. Moreover, as the degree of the polynomials to be reproduced increases, the smoothness of the basis functions also improves. As shown in Theorem \ref{thm_convergence_polynomial_reproduction_exponential_decay}, the smoothness of the weight functions does not influence the convergence rate. However, smooth approximants can be advantageous in many applications.
\\
\\
To generate Figures \ref{fig_exp_2_basis_function_fast_decay} and Figure \ref{fig_exp_basis_function_fast_decay} , we solved the linear system in equation \eqref{eq_del_lambda_k_x_exp_fast_decay} using Chebyshev polynomials of the first kind \cite{Rivlin_1990} as the polynomial basis.
\\

\subsection{Numerical Convergence}

We present a series of numerical experiments aimed at verifying numerically Theorem \ref{thm_convergence_polynomial_reproduction_exponential_decay}. For a comprehensive study of the proposed methods in the one-dimensional case, we refer the reader to \cite{cappellazzo_2022}.
In principle, the proposed approximation methods can be applied to data sets of arbitrary size. The moving least squares scheme \eqref{eq_moving_least_squares_fast_decay} have been implemented in two dimensions. We approximate Franke's bivariate test function \cite{Franke_1979} on a set of equispaced (quasi-uniform) nodes in the unit square $[0,1]^2$.
The approximant is computed via the linear system in equation \eqref{eq_del_lambda_k_x_exp_fast_decay} using the standard monomial basis.  
We set $\delta=30h_{X,\Omega}$, where the parameter was selected through a trial-and-error procedure. Compared to the one-dimensional case, our numerical experiments indicate that the choice of the parameter $\delta$ plays a more critical role in two dimensions.
\\
\\
As the weight function, we used
\begin{equation*}
    w(x,y)=e^{-\left(\frac{|x-y|}{\delta}\right)^2}.
\end{equation*}
In the following graphs, the blue line serves as a reference to verify the expected slope of the approximation error, in accordance with Theorem \ref{thm_convergence_polynomial_reproduction_exponential_decay}.

\begin{figure}[H]
	\begin{subfigure}{0.45\textwidth}
		\includegraphics[width=\linewidth]{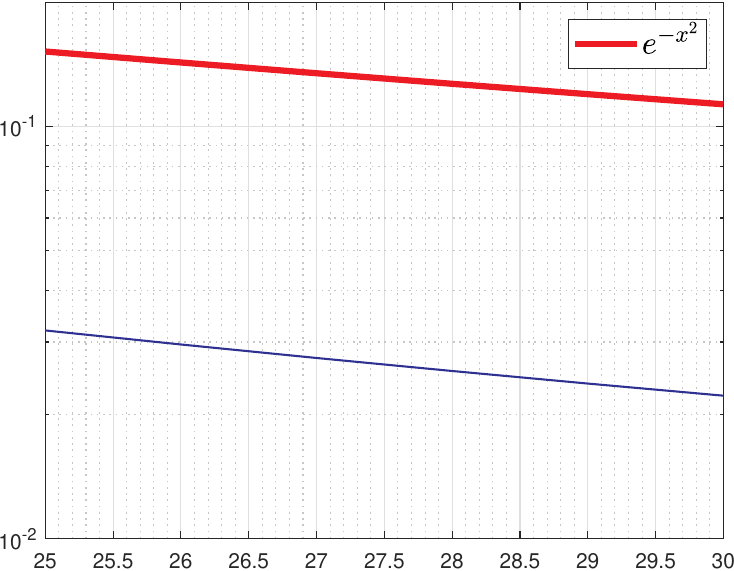}
		%		\label{fig:1}
  \end{subfigure}\hfil % <-- added
  \begin{subfigure}{0.45\textwidth}
		\includegraphics[width=\linewidth]{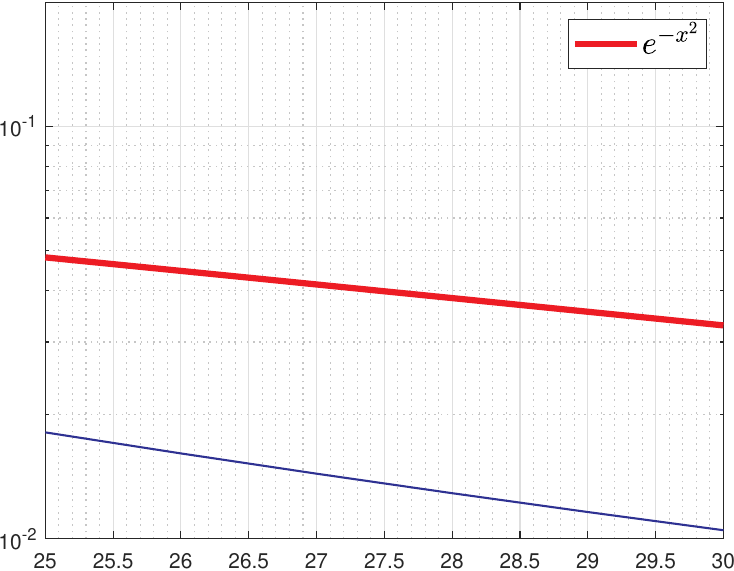}
		%		\label{fig:1}
  \end{subfigure}\hfil % <-- added
	
\caption{Convergence rate of the approximation error $\|f-z_{f,X}\|_{L^{\infty}([-1,1])}$. The x-axis represents the number of equispaced nodes used to construct the approximant. The approximation method reproduces polynomials of degree $1$ (left) and $2$ (right) exactly.} \label{franke_errors_slode_exp_2}
\end{figure}

\begin{table}[H]
\begin{center}
\begin{tabular}{||c||c|c|c|c|c|c||c||}
\hline
Nodes & 676 & 729 & 784 & 841 & 900 & 961 & Degree \\ \hline
Left & 1.52e-1 & 1.43e-1 & 1.35e-1 & 1.27e-1 & 1.20e-1 & 1.13e-1 & 1 \\ \hline
Right & 4.81e-2 & 4.46e-2 & 4.14e-2 & 3.83e-2 & 3.55e-2 & 3.29e-2 & 2 \\ \hline
%$\| \cdot \|_1 - e^{-x^2}$ & 1.68e-2 & 1.55e-2 & 1.48e-2 & 1.35e-2 & 1.29e-2 & 9.66e-3 & 1 \\ \hline
\end{tabular}
\end{center}
\caption{Approximation errors $\|f-z_{f,X}\|_{L^{\infty}([-1,1])}$ relative to Figure \ref{franke_errors_slode_exp_2}.}
\end{table}

\noindent The previous numerical experiment confirms the statement of Theorem \ref{thm_convergence_polynomial_reproduction_exponential_decay} also in the multivariate setting.

\subsection{Basis functions with algebraic decay}
\label{subsection_algebraic_decay_examples}

As shown in equation \eqref{eq_shepard_method_general}, the Shepard method constitutes a particular case of the proposed approximation schemes. It is therefore of interest to examine basis functions with algebraic decay.
\\
\\
In this section, we consider weight functions satisfying
\begin{equation}
    \varphi(x) \leq \frac{1}{x^k} \quad \text{ for } \, x \in \mathbb{R}_{>0} \quad \textrm{ and } \quad k \in \mathbb{R}_{>0}.
    \label{eq_algebraic_decay_def}
\end{equation}
In this setting, $\varphi$ may exhibit a singularity at the origin and may fail to satisfy the ratio test (cf. Definition \ref{def_fast_decay_pol_rep}). As shown in Theorem \ref{thm_convergence_polynomial_reproduction_exponential_decay}, the convergence of the fast-decaying polynomial reproduction framework relies on inequality \eqref{eq_stability_quasi_interpoaltion_fast_decay}.
\\
\\
From inequalities 
\begin{equation*}
    |u_j(x)| \leq C\varphi\left(\frac{\|x-x_j\|_2}{q_X}\right) \leq C \frac{q_X^k}{\|x-x_j\|_2^k} \; \text{ for all } \; x \in \Omega \; \text{ and } \; j=1,\dots,N
\end{equation*}
we can derive the inequality \eqref{eq_stability_quasi_interpoaltion_fast_decay} with a (possibly different) constant $K$. The proof follows the same steps as that of Theorem \ref{thm_stability_quasi_interpoaltion_fast_decay}, with the only modification being that algebraic decay is used when analyzing the node distribution in the set $E_0$.
\\
The choice of the exponent $k$ in \eqref{eq_algebraic_decay_def} is critical for both the convergence and stability of the method. If the method reproduces polynomials up to degree $m$ in $\mathbb{R}^d$, then the following conditions must hold: $d+m-\frac{k}{2}<-1$ for the moving least squares method.
%, and $d+m-k<-1$ for the $1$-norm approximation. 
This condition was derived from equation \eqref{eq_general_phi_weight_mls} %and \eqref{eq_general_phi_weight_1_norm} 
together with the convergence of the generalized harmonic series. This constraint on the exponent $k$ is consistent with those found in \cite{LightCheney1992}. Indeed, there exists $\lambda \in ]0,1[$ such that, in a neighborhood of $+\infty$ that does not contain $0$,  say for instance $]1, +\infty[$
\begin{equation*}
    \sup_{x\in ]1,+ \infty[}  \left\{|\varphi(x)| (1+ x)^{d+m+1+\lambda} \right\}  < +\infty.
\end{equation*}
With regard to algebraic decay, the proposed schemes follow the construction introduced in \cite{LightCheney1992} while extending it to allow for possible divergence at the origin. A further advantage of our framework is that it enables us to relate the regularity of the approximant to the smoothness of the basis functions, thereby extending the range of applicability beyond that of continuous functions with algebraic decay.
\\
\\
As  test function we used
\begin{equation*}
    f(x) = \sin(\pi x).
\end{equation*}

\begin{figure}[H]
\begin{center}
	\begin{subfigure}{0.45\textwidth}
		\includegraphics[width=\linewidth]{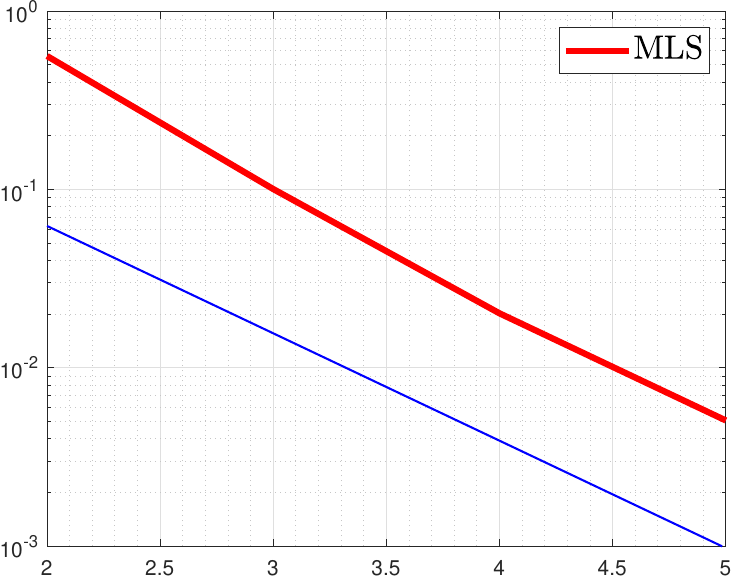}
				\caption{$k=6.2$}
  \end{subfigure}\hfil % <-- added
\end{center}
	
\caption{Convergence rate of the approximation error $\|f-z_{f,X}\|_{L^{\infty}([-1,1])}$. The x-axis represents the number of equispaced nodes used to construct the approximant. The approximation method reproduces polynomials of degree $1$ exactly.  In this experiment, the scaling parameter is set to $\delta = 5h_{X,\Omega}$.} \label{sin_algebraic_decay_example_001}
\end{figure}

\begin{table}[H]
\begin{center}
\begin{tabular}{||c||c|c|c|c||c||}
\hline
Nodes & 8 & 16 & 32 & 64 & Degree \\ \hline
MLS & 5.61e-01 & 1.01e-01 & 2.02e-02 & 5.07e-03 & 1 \\ \hline
%$\|\cdot\|_1$ & 6.82e-02 & 1.52e-02 & 4.41e-03 & 9.31e-04 & 1 \\ \hline
\end{tabular}
\end{center}
\caption{Approximation errors $\|f-z_{f,X}\|_{L^{\infty}([-1,1])}$ relative to Figure \ref{sin_algebraic_decay_example_001}.}
\end{table}
\noindent Following these considerations, it is straightforward to observe that the fast-decaying polynomial reproduction framework can be applied to other kernels without compact support. For instance, the inverse multiquadric radial basis function \cite{hardy_71} can be employed as a weight function:
\begin{equation*}
    \varphi(x)=(c^2+x^2)^{\beta}, \qquad \beta<0.
\end{equation*}
We can recognize the previous construction because
\begin{equation*}
    (c^2+x^2)^{\beta} \leq (Mx^2)^{\beta},
\end{equation*}
whenever $M \in ]0,1[$.

\subsection{Stability}

The objective of this section is to numerically investigate the Lebesgue constant and the stability of the approximation methods studied in this work. Theorem \ref{thm_stability_quasi_interpoaltion_fast_decay} enables us to estimate the Lebesgue constant independently of the parameters appearing in equations \eqref{eq_delta_h_X_relation} and \eqref{eq_delta_q_X_relation} (quasi-uniform node distribution). Throughout this section, we set $\delta=5h_{X,\Omega}$.
\\
\\
As the weight function, we used
\begin{equation*}
    w(x,y)=e^{-\left(\frac{|x-y|}{\delta}\right)^2}.
\end{equation*}

\begin{figure}[H]
\begin{center}
	\begin{subfigure}{0.45\textwidth}
		\includegraphics[width=\linewidth]{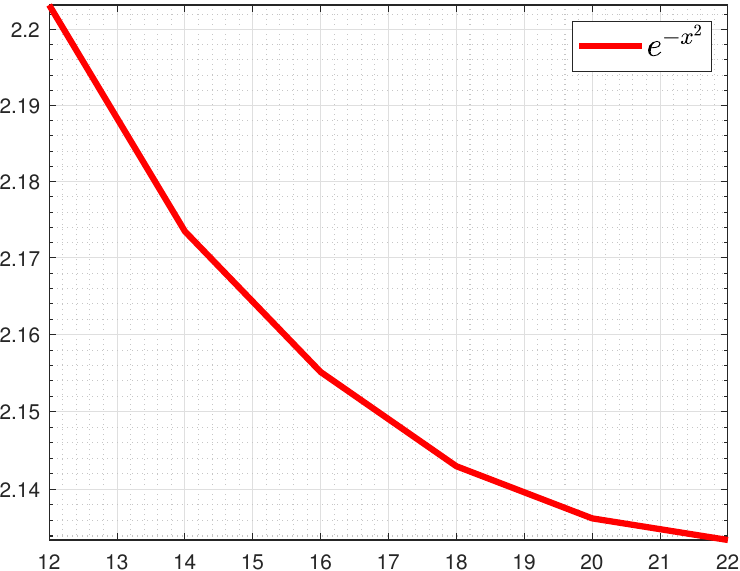}
		%		\label{fig:1}
  \end{subfigure}\hfil % <-- added
\end{center}
\caption{Lebesgue constant of the moving least squares method in equation \eqref{eq_moving_least_squares_fast_decay}. %and the $1$-norm minimization method in equation \eqref{approximant_linear_problem_optimization_form}. 
The x-axis represents the number of equispaced nodes used to construct the basis functions. The approximation method reproduce polynomials of degree $2$ exactly.} \label{lebesgue_constant_example_exp_2_001}
\end{figure}

\begin{table}[H]
\begin{center}
\begin{tabular}{||c||c|c|c|c|c|c||c||}
\hline
Nodes & 169 & 225 & 289 & 361 & 441 & 529 & Degree \\ \hline
MLS & 2.203275 & 2.173473 & 2.155158 & 2.142972 & 2.136264 & 2.133485 & 2 \\ \hline
%$\| \cdot \|_1 - e^{-x^2}$ & 1.017143 & 1.016552 & 1.016274 & 1.016000 & 1.016361 & 1.017165 & 2 \\ \hline
\end{tabular}
\end{center}
\caption{Lebesgue constant of the moving least squares method in equation \eqref{eq_moving_least_squares_fast_decay} relative to Figure \ref{lebesgue_constant_example_exp_2_001}.}
\end{table}

\section{Conclusion}

The rescaled localized radial basis function (RL-RBF) method, proposed in \cite{deparis_2014}, provides an efficient numerical scheme for interpolating functions on scattered nodes. A proof of convergence for the RL-RBF method was later given in \cite{demarchi2020}, although it relies on a conjecture in the quasi-uniform setting (cf. eq. \eqref{RL_RBF_approximant_definition_001}, a lower bound for the RBF approximation of the constant function $1$ can be derived in the \textit{stationary setting}, regardless of the distribution of the nodes).
\\
\\
Our work extends this approach by generalizing the RL-RBF method through an increase in the dimension of the polynomial space that is reproduced exactly. The objective is to construct a convergent approximation method whose convergence rate is of order $\mathcal{O}(h_{X,\Omega}^{m+1})$ whenever all polynomials of degree at most $m$ are reproduced exactly. To this end, we modify the notion of local polynomial reproduction by replacing the requirement of compact support with that of fast decay of the basis functions. The RL-RBF method fits naturally into this new framework. %, also because its cardinal functions do not have compact support. 
In the quasi-uniform setting, the resulting approximation scheme is both convergent and stable, with the Lebesgue constant depending on the ambient dimension, the degree of the polynomial space being reproduced, and the decay rate of the basis functions.
\\
\\
At this stage of the analysis, smoothness does not play a fundamental role. In addition to the RL-RBF method, we propose a further approximation scheme that yields smooth quasi-interpolants. This method approximates the value of an unknown function via the moving least squares technique with a Gaussian weight function. The smoothness of the resulting approximant is inherited from that of the weight function. Since the solution of a moving least squares problem is obtained by solving a quadratic optimization problem, the computational cost of the method is linear in the number of nodes. Numerical experiments confirm the theoretical results on convergence and stability, as well as the smoothness of the basis functions. An unexpected observation arising from the numerical tests is that, even when the weight function is only continuous, increasing the degree of the polynomial space to be reproduced renders the basis functions numerically smooth. 
The methods have also been tested in the multivariate setting. We have further experimented with degenerate basis functions exhibiting algebraic decay and a singularity at the origin. The numerical results support the validity of the framework, although the choice of the decay exponent proves critical for achieving both stability and convergence.
\subsection{Future Work}
In the preceding analysis and numerical tests, we considered functions with global support, resulting in potentially dense matrices even when the polynomial space to be reproduced is not too large. Future work is directed toward mitigating this difficulty: we can replace the quadratic optimization problem arising from moving least squares (eq. \eqref{eq_def_minimization_mls_2_exp_fast_decay}) with a linear program over a polyhedron. A vertex solution of the linear program allows us to control the number of nonzero basis functions at each point in the domain. Since the weight functions are designed to localize the optimization problem, we expect the weight corresponding to a given node to be large when the evaluation point is far from that node. This observation motivates the use of \textit{column generation techniques} to reduce the dimensionality of the problem, since the number of nonzero basis functions is uniformly bounded with respect to the dimension of the reproduced polynomial space.
\\
\\
We plan to extend the applicability of fast-decaying polynomial reproduction schemes to more general point distributions by employing the Fake Nodes Approach \cite{demarchi_2021}. Another promising direction concerns the optimality of the Lebesgue constant appearing in equation \eqref{def_constant_K_stability}. Indeed, by appropriately tuning the parameters $C$ and $\varphi$, it may be possible to improve stability and achieve optimal control over the decay of the basis functions.
\vskip 0.1in 

\noindent {\bf Acknowledgments}. 

\noindent This work was initiated during an Erasmus traineeship at the University of Bayreuth under the supervision of Professor Holger Wendland. The authors are grateful to Professor Martin Buhmann of the University of Giessen for many instructive discussions. This research was carried out as part of the UMI topic group ``Teoria dell'Approssimazione e Applicazioni''. The authors are also members of the INdAM-GNCS research group and the SIMAI ANA\&A activity group.

\printbibliography %Prints bibliography
%\bibliographystyle{alpha}
%\biblography{reference\_FDPR}

\end{document}